\documentclass [12pt]{article}
\usepackage{amsmath,amssymb}

\begin{document}

   \title {Analogues of the adjoint matrix for generalized inverses
    and corresponding Cramer rules.}

\author {Kyrchei I. I.\footnote{Pidstrygach Institute
for Applied Problems of Mechanics and Mathematics, str.Naukova 3b,
Lviv, Ukraine, 79005, kyrchei@lms.lviv.ua}}
\date{}
\maketitle
\begin{abstract}
 In this article, we introduce determinantal
representations
 of the Moore - Penrose inverse and the Drazin inverse which are based on analogues  of the
classical adjoint matrix. Using the obtained analogues  of the
adjoint matrix, we get  Cramer rules for the least squares
solution and for the Drazin inverse solution of singular linear
systems. Finally, determinantal expressions for ${\rm {\bf A}}^{
+} {\rm {\bf A}}$, ${\rm {\bf A}} {\rm {\bf A}}^{ +}$, and ${\rm
{\bf A}}^{ D} {\rm {\bf A}}$ are presented.
\end{abstract}
%\begin
\textit{Keywords}: Moore-Penrose inverse; Drazin inverse; system
of linear equations; least squares solution; Cramer rule

\noindent \textit{AMS classification}: 15A09, 15A57
%\end{keyword}
%\end{frontmatter}

\section{Introduction}
   Determinantal representation of the
Moore - Penrose inverse was studied in \cite{ba, ben1,ga,mo,st1}.
The main result  consists in the following theorem.
\newtheorem{theorem}{Theorem}[section]
\begin{theorem}\label{kyrc1} The Moore - Penrose inverse
 ${\rm {\bf A}}^{+}=(a_{ij}^{+})_{n\times m}$
of ${\rm {\bf A}}\in {\mathbb C}_{r}^{m\times n} $ has the
following determinantal representation
\[
a_{ij}^{+} = {\frac{{{\sum\limits_{\left( {\alpha,\,\beta} \right)
\in N_{r}  {\left\{ {j,\,i} \right\}}} {{\left| {\left( {{\rm {\bf
A}}^{\ast}} \right)_{\alpha} ^{\beta} } \right|}{\frac{{\partial}
}{{\partial a_{j\,i} }}}{\left| {{\rm {\bf A}}_{\beta} ^{\alpha} }
\right|}} }}}{{{\sum\limits_{\left( {\gamma,\,\delta}  \right) \in
N_{r} } {{\left| {\left( {{\rm {\bf A}}^{\ast}} \right)_{\gamma}
^{\delta} } \right|}\,{\left| {{\rm {\bf A}}_{\delta} ^{\gamma} }
\right|}}} }}}, \,\,1 \leq i, j \leq n.
\]
 \end{theorem}
 Stanimirovic' \cite{st2} introduced a determinantal representation of the Drazin inverse
by the following theorem.
\begin{theorem}
The Drazin inverse ${\rm {\bf A}}^{D}=\left(a_{ij}^{D}\right)$ of
an arbitrary matrix ${\rm {\bf A}}\in {\mathbb C}^{n\times n} $
with $Ind{\rm {\bf A}}=k $ possesses the following determinantal
representation
\[
a_{ij}^{D} = {\frac{{{\sum\limits_{\left( {\alpha,\,\beta} \right)
\in N_{r_{k}}  {\left\{ {j,\,i} \right\}}} {{\left| {\left( {{\rm
{\bf A}}^{s}} \right)_{\alpha} ^{\beta} }
\right|}{\frac{{\partial} }{{\partial a_{j\,i} }}}{\left| {{\rm
{\bf A}}_{\beta} ^{\alpha} }  \right|}} }}}{{{\sum\limits_{\left(
{\gamma,\,\delta}  \right) \in N_{r_{k}} } {{\left| {\left( {{\rm
{\bf A}}^{s}} \right)_{\gamma} ^{\delta} } \right|}\,{\left| {{\rm
{\bf A}}_{\delta} ^{\gamma} }  \right|}}} }}}, \,\,1 \leq i, j
\leq n;
\]
where $s \ge k$ and $r_{k} = rank{\rm {\bf A}}^{s}$. \end{theorem}
These determinantal representations  of generalized  inverses are
based on corresponding full-rank representations.

 We use the
following notations from \cite{ba,st1}. Let ${\mathbb C}^{m\times
n} $ be the set of $m$ by $n$ matrices with complex entries,
${\mathbb C}^{m\times n}_{r} $ be the subset of ${\mathbb
C}^{m\times n} $  in which every matrix has rank $r$. ${\rm \bf
I}_{m}$ denotes the identity matrix of order $m$, and
$\|.\|=\|.\|_{2}$ is the Euclidean vector norm.  Let $\alpha : =
\left\{ {\alpha _{1} ,\ldots ,\alpha _{k}}  \right\} \subseteq
{\left\{ {1,\ldots ,m} \right\}}$ and $\beta : = \left\{ {\beta
_{1} ,\ldots ,\beta _{k}} \right\} \subseteq {\left\{ {1,\ldots
,n} \right\}}$ be subsets of the order $1 \le k \le \min {\left\{
{m,n} \right\}}$. Then ${\left| {{\rm {\bf A}}_{\beta} ^{\alpha} }
\right|}$ denotes the minor of  ${\rm {\bf A}}\in {\mathbb
C}^{m\times n} $ determined by the rows indexed by $\alpha$ and
the columns indexed by $\beta$. Clearly, ${\left| {{\rm {\bf
A}}_{\alpha} ^{\alpha} } \right|}$ denotes a principal minor
determined by the rows and columns indexed by $\alpha$.  The
cofactor of $a_{ij} $ in  ${\rm {\bf A}}\in {\mathbb C}^{n\times
n} $ is denoted by ${\frac{{\partial }}{{\partial a_{ij}}
}}{\left| {{\rm {\bf A}}} \right|}$. For $1 \leq k\leq n$, denote
by $\textsl{L}_{ k, n}: = {\left\{ {\,\alpha :\alpha = \left(
{\alpha _{1} ,\ldots ,\alpha _{k}} \right),\,{\kern 1pt} 1 \le
\alpha _{1} \le \ldots \le \alpha _{k} \le n} \right\}}$ the
collection of strictly increasing sequences of $k$ integers chosen
from $\left\{ {1,\ldots ,n} \right\}$. Let $N_{k} : =
\textsl{L}_{k, m} \times \textsl{L}_{ k, n} $. For fixed $\alpha
\in \textsl{L}_{ p, m}$, $\beta \in \textsl{L}_{ p, n}$, $1\leq
p\leq k$, let
\[\begin{array}{c}
   I_{k,\,m} \left( {\alpha} \right): = {\left\{ {I:\,I \in
\textsl{L}_{ k,\, m} ,I \supseteq \alpha} \right\}},\\
  J_{k,\,n}
\left( {\beta} \right): = {\left\{ {J:\,J \in \textsl{L}_{ k,\, n}
,J \supseteq \beta} \right\}},\\
N_{k} \left( {\alpha ,\beta} \right): = I_{k,\,m} \left( {\alpha}
\right)\times J_{k,\,n} \left( {\beta}  \right)
\end{array}
\]
For case $i \in \alpha $ and $j \in \beta$, we denote
\[\begin{array}{c}
I_{k,m} {\left\{ {i} \right\}}: = {\left\{ {\alpha:\,\alpha \in
L_{k, m},  i \in \alpha} \right\}},
 J_{k,\,n} {\left\{ {j}
\right\}}: =  {\left\{ {\beta:\,\beta\in  L_{k, n}, j \in  \beta}
\right\}},\\
N_{k} {\left\{ {i,j} \right\}}: = I_{k,\, m} {\left\{ {i}
\right\}}\times J_{k,\, n} {\left\{ {j} \right\}}.
\end{array}
\]

 In this paper
we introduce  determinantal representations of the Moore - Penrose
inverse and of the Drazin inverse  based on corresponding limit
representations. The obtained determinantal representations can be
considered as founded on some analogues of the classical adjoint
matrix.  The corresponding Cramer rules for the complex system of
linear equations with a rectangular or singular coefficient matrix
follow  from these analogues.

\section{ Analogues of the classical adjoint matrix for the Moore
- Penrose inverse}
 We shall use the following well-known facts
(see, for example, \cite{ho}).
\newtheorem{definition}{Definition}[section]
\begin{definition}
The matrix ${\rm {\bf A}}^{ +} \in {\mathbb C}^{n\times m}$ is
called the Moore - Penrose inverse of an arbitrary ${\rm {\bf
A}}\in {\mathbb C}^{m\times n}$ if it satisfies the equations
\[{\rm {\bf A}}{\rm {\bf A}}^{ +} {\rm {\bf A}} = {\rm {\bf
A}};\,\,{\rm {\bf A}}^{ +}{\rm {\bf A}}{\rm {\bf A}}^{ +}  = {\rm
{\bf A}}^{ +};\,\, ({\rm {\bf A}}{\rm {\bf A}}^{ + })^{\ast}={\rm
{\bf A}}{\rm {\bf A}}^{ + };\,\,({\rm {\bf A}}^{ +} {\rm {\bf
A}})^{\ast}={\rm {\bf A}}^{ +} {\rm {\bf A}}.\]
\end{definition}
The superscript $\ast$ denotes conjugate transpose matrix.
\newtheorem{lemma}{Lemma}[section]
 \begin{lemma} \cite{ho} There exists a unique
Moore - Penrose inverse ${\rm {\bf A}}^{ + }$ of ${\rm {\bf A}}\in
{\mathbb C}^{m\times n} $.
\end{lemma}
 \begin{lemma}\cite{ho}\label{kyrc2}
If ${\rm {\bf A}}\in {\mathbb C}^{m\times n} $, then
\[{\rm {\bf
A}}^{ +} = {\mathop {\lim }\limits_{\lambda \to 0}} {\rm {\bf
A}}^{ *} \left( {{\rm {\bf A}}{\rm {\bf A}}^{ *}  + \lambda {\rm
{\bf I}}} \right)^{ - 1} = {\mathop {\lim }\limits_{\lambda\to 0}}
\left( {{\rm {\bf A}}^{ *} {\rm {\bf A}} + \lambda {\rm {\bf I}}}
\right)^{ - 1}{\rm {\bf A}}^{ *}, \]
 where $\lambda \in {\mathbb
R} _{ +}  $, and ${\mathbb R} _{ +} $ is the set of  positive real
numbers.
\end{lemma}
 \begin{lemma} \cite{ho}\label{kyrc3} If
${\rm {\bf A}}\in {\mathbb C}^{m\times n} $, then the following
statements are true.
 \begin{itemize}
\item [ i)] If $\rm{rank}\,{\rm {\bf A}} = n$, then ${\rm {\bf A}}^{ +}
= \left( {{\rm {\bf A}}^{ *} {\rm {\bf A}}} \right)^{ - 1}{\rm
{\bf A}}^{ * }$ .
\item [ ii)] If $\rm{rank}\,{\rm {\bf A}} =
m$, then ${\rm {\bf A}}^{ +}  = {\rm {\bf A}}^{ * }\left( {{\rm
{\bf A}}{\rm {\bf A}}^{ *} } \right)^{ - 1}.$
\item [ iii)] If $\rm{rank}\,{\rm {\bf A}} = n = m$, then ${\rm {\bf
A}}^{ +}  = {\rm {\bf A}}^{ - 1}$ .
\end{itemize}
\end{lemma}
\begin{theorem} \cite{ho}\label{kyrc4} Let $d_{r}$ be the sum of  principal minors
of order $r$ of  ${\rm {\bf A}}\in {\mathbb C}^{n\times n}$. Then
its characteristic polynomial $ p_{{\rm {\bf A}}}\left( {t}
\right)$ can be expressed as  $ p_{{\rm {\bf A}}}\left( {t}
\right) = \det \left( { t{\rm {\bf I}} - {\rm {\bf A}}} \right) =
t^{n} - d_{1} t^{n - 1} + d_{2} t^{n - 2} - \ldots + \left( { - 1}
\right)^{n}d_{n}. $
\end{theorem}
Denote by ${\rm {\bf a}}_{.j} $ and ${\rm {\bf a}}_{i.} $ the
$j$th column  and the $i$th row of  ${\rm {\bf A}}\in {\mathbb
C}^{m\times n} $ respectively. In the same way, denote by ${\rm
{\bf a}}^{\ast}_{.j} $ and ${\rm {\bf a}}^{\ast}_{i.} $ the $j$th
column  and the $i$th row of Hermitian adjoint matrix ${\rm {\bf
A}}^{\ast}$. Let ${\rm {\bf A}}_{.j} \left( {{\rm {\bf b}}}
\right)$ denote the matrix obtained from ${\rm {\bf A}}$ by
replacing its $j$th column with some  vector ${\rm {\bf b}}$, and
let ${\rm {\bf A}}_{i.} \left( {{\rm {\bf b}}} \right)$ denote the
matrix obtained from ${\rm {\bf A}}$ by replacing its $i$th row
with ${\rm {\bf b}}$.

\newcommand{\rank}{\mathop{\rm rank}\nolimits}
\begin{lemma} \label{kyrc5} If ${\rm {\bf A}}\in {\mathbb C}^{m\times n}_{r}
$, then $
 \rank\,\left( {{\rm {\bf A}}^{ *} {\rm {\bf A}}}
\right)_{.\,i} \left( {{\rm {\bf a}}_{.j}^{ *} }  \right) \le r.
$
\end{lemma}
 {\textit{Proof}}.
Let ${\rm {\bf P}}_{i\,k} \left( {-a_{j\,k}}  \right)\in {\mathbb
C}^{n\times n} $, $(k \ne i )$, be the matrix with $-a_{j\,k} $ in
the $(i, k)$ entry, 1 in all diagonal entries, and 0 in others. It
is the  matrix of an elementary transformation. It follows that
\[
\left( {{\rm {\bf A}}^{ *} {\rm {\bf A}}} \right)_{.\,i} \left(
{{\rm {\bf a}}_{.\,j}^{ *} }  \right) \cdot {\prod\limits_{k \ne
i} {{\rm {\bf P}}_{i\,k} \left( {-a_{j\,k}}  \right) = {\mathop
{\left( {{\begin{array}{*{20}c}
 {{\sum\limits_{k \ne j} {a_{1k}^{ *}  a_{k1}} } } \hfill & {\ldots}  \hfill
& {a_{1j}^{ *} }  \hfill & {\ldots}  \hfill & {{\sum\limits_{k \ne
j } {a_{1k}^{ *}  a_{kn}}}} \hfill \\
 {\ldots}  \hfill & {\ldots}  \hfill & {\ldots}  \hfill & {\ldots}  \hfill &
{\ldots}  \hfill \\
 {{\sum\limits_{k \ne j} {a_{nk}^{ *}  a_{k1}} } } \hfill & {\ldots}  \hfill
& {a_{nj}^{ *} }  \hfill & {\ldots}  \hfill & {{\sum\limits_{k \ne
j } {a_{nk}^{ *}  a_{kn}}}} \hfill \\
\end{array}} }
\right)}\limits_{i-th}}}}.
\]
 The obtained above matrix  has the following factorization.
\[
{\mathop {\left( {{\begin{array}{*{20}c}
 {{\sum\limits_{k \ne j} {a_{1k}^{ *}  a_{k1}} } } \hfill & {\ldots}  \hfill
& {a_{1j}^{ *} }  \hfill & {\ldots}  \hfill & {{\sum\limits_{k \ne
j } {a_{1k}^{ *}  a_{kn}} } } \hfill \\
 {\ldots}  \hfill & {\ldots}  \hfill & {\ldots}  \hfill & {\ldots}  \hfill &
{\ldots}  \hfill \\
 {{\sum\limits_{k \ne j} {a_{nk}^{ *}  a_{k1}} } } \hfill & {\ldots}  \hfill
& {a_{nj}^{ *} }  \hfill & {\ldots}  \hfill & {{\sum\limits_{k \ne
j } {a_{nk}^{ *}  a_{kn}} } } \hfill \\
\end{array}} } \right)}\limits_{i-th}}  =
\]
\[
 = \left( {{\begin{array}{*{20}c}
 {a_{11}^{ *} }  \hfill & {a_{12}^{ *} }  \hfill & {\ldots}  \hfill &
{a_{1m}^{ *} }  \hfill \\
 {a_{21}^{ *} }  \hfill & {a_{22}^{ *} }  \hfill & {\ldots}  \hfill &
{a_{2m}^{ *} }  \hfill \\
 {\ldots}  \hfill & {\ldots}  \hfill & {\ldots}  \hfill & {\ldots}  \hfill
\\
 {a_{n1}^{ *} }  \hfill & {a_{n2}^{ *} }  \hfill & {\ldots}  \hfill &
{a_{nm}^{ *} }  \hfill \\
\end{array}} } \right){\mathop {\left( {{\begin{array}{*{20}c}
 {a_{11}}  \hfill & {\ldots}  \hfill & {0} \hfill & {\ldots}  \hfill &
{a_{n1}}  \hfill \\
 {\ldots}  \hfill & {\ldots}  \hfill & {\ldots}  \hfill & {\ldots}  \hfill &
{\ldots}  \hfill \\
 {0} \hfill & {\ldots}  \hfill & {1} \hfill & {\ldots}  \hfill & {0} \hfill
\\
 {\ldots}  \hfill & {\ldots}  \hfill & {\ldots}  \hfill & {\ldots}  \hfill &
{\ldots}  \hfill \\
 {a_{m1}}  \hfill & {\ldots}  \hfill & {0} \hfill & {\ldots}  \hfill &
{a_{mn}}  \hfill \\
\end{array}} } \right)}\limits_{i-th}} j-th.
\]
 Denote by ${\rm {\bf \tilde {A}}}: = {\mathop
{\left( {{\begin{array}{*{20}c}
 {a_{11}}  \hfill & {\ldots}  \hfill & {0} \hfill & {\ldots}  \hfill &
{a_{1n}}  \hfill \\
 {\ldots}  \hfill & {\ldots}  \hfill & {\ldots}  \hfill & {\ldots}  \hfill &
{\ldots}  \hfill \\
 {0} \hfill & {\ldots}  \hfill & {1} \hfill & {\ldots}  \hfill & {0} \hfill
\\
 {\ldots}  \hfill & {\ldots}  \hfill & {\ldots}  \hfill & {\ldots}  \hfill &
{\ldots}  \hfill \\
 {a_{m1}}  \hfill & {\ldots}  \hfill & {0} \hfill & {\ldots}  \hfill &
{a_{mn}}  \hfill \\
\end{array}} } \right)}\limits_{i-th}} j-th$. The matrix
${\rm {\bf \tilde {A}}}$ is obtained from ${\rm {\bf A}}$ by
replacing all entries of the $j$th row  and of the $i$th column
with zeroes except that the $(j, i)$ entry equals 1. Elementary
transformations of a matrix do not change its rank. It follows
that $\rank\left( {{\rm {\bf A}}^{ *} {\rm {\bf A}}}
\right)_{.\,i} \left( {{\rm {\bf a}}_{.j}^{ *} } \right) \le \min
{\left\{ {\rank{\rm {\bf A}}^{ * },\rank{\rm {\bf \tilde {A}}}}
\right\}}$. Since $\rank{\rm {\bf \tilde {A}}} \ge \rank\,{\rm
{\bf A}} = \rank{\rm {\bf A}}^{ *} $ and $\rank{\rm {\bf A}}^{
* }{\rm {\bf A}} = \rank{\rm {\bf A}}$ the proof is completed.$\blacksquare$

The following lemma is  proved in the same way.
\begin{lemma} If ${\rm {\bf A}}\in {\mathbb C}^{m\times n}_{r} $,
then $ \rank\left( {{\rm {\bf A}}{\rm {\bf A}}^{ *} }
\right)_{i\,.} \left( {{\rm {\bf a}}_{j\,.}^{ *} }  \right) \le
r.$
\end{lemma}
\begin{theorem}\label{kyrc6}
The Moore-Penrose inverse ${\rm {\bf A}}^{+}$ of ${\rm {\bf A}}\in
{\mathbb C}^{m\times n}_{r} $ can be represented  as follows
\begin{equation}\label{eq1}
\begin{array}{l}
  {\rm {\bf A}}^{
+}=\left({\frac{l_{ij}}{{d_{r} ({\rm {\bf A}}^{ *} {\rm {\bf
A}})}}}
 \right)_{n\times m},\,
\mbox{where}
\\
l_{ij} = {\sum\limits_{\beta \in J_{r,\, n} {\left\{ {i}
\right\}}} {{\left| {\left( {({\rm {\bf A}}^{ *} {\rm {\bf
A}})_{.\,i} ({\rm {\bf a}}_{.j}^{ *}  )} \right)_{\beta} ^{\beta}
}  \right|}}},\,\,\,
  d_{r} ({\rm {\bf A}}^{ *} {\rm {\bf A}}) =
{\sum\limits_{\beta \in J_{r, \,n}}  {{\left| {\left( {{\rm {\bf
A}}^{ *} {\rm {\bf A}}} \right)_{\beta} ^{\beta} } \right|}}},
\end{array}
\end{equation}
or \begin{equation}\label{eq2}
 \begin{array}{l}
   {\rm {\bf A}}^{ +}  =\left({\frac{{r_{ij}}}{{d_{r}
 ({\rm {\bf A}}{\rm {\bf A}}^{ *} )}}} \right)_{n\times
 m},\,
 \mbox{where} \\
   r_{ij} =
{\sum\limits_{\alpha \in I_{r,\, m} {\left\{ {j} \right\}}}
{{\left| {\left( {({\rm {\bf A}}{\rm {\bf A}}^{ * })_{j\,.} ({\rm
{\bf a}}_{i\,.}^{ *}  )} \right)_{\alpha} ^{\alpha} }
\right|}}},\,\,\,  d_{r} ({\rm {\bf A}}{\rm {\bf A}}^{ *} ) =
{\sum\limits_{\alpha \in I_{r, \,m} } {{\left| {\left( {{\rm {\bf
A}}{\rm {\bf A}}^{ *} } \right)_{\alpha }^{\alpha} } \right|}}} .
 \end{array}
\end{equation}

\end{theorem}
{\textit{Proof}}. At first we shall obtain the representation
(\ref{eq1}). If $\lambda \in {\mathbb R} _{ +}  $, then the matrix
$\left( {\lambda {\rm {\bf I}} + {\rm {\bf A}}^{ *} {\rm {\bf A}}}
\right)\in {\mathbb C}^{n\times n} $ is Hermitian and $\rank
\left( {\lambda {\rm {\bf I}} + {\rm {\bf A}}^{ *} {\rm {\bf A}}}
\right)=n$.
 Hence,   there exists its  inverse
\[
\left( {\lambda {\rm {\bf I}} + {\rm {\bf A}}^{ *} {\rm {\bf A}}}
\right)^{ - 1} = {\frac{{1}}{{\det \left( {\lambda {\rm {\bf I}} +
{\rm {\bf A}}^{ * }{\rm {\bf A}}} \right)}}}\left(
{{\begin{array}{*{20}c}
 {L_{11}}  \hfill & {L_{21}}  \hfill & {\ldots}  \hfill & {L_{n\,1}}  \hfill
\\
 {L_{12}}  \hfill & {L_{22}}  \hfill & {\ldots}  \hfill & {L_{n\,2}}  \hfill
\\
 {\ldots}  \hfill & {\ldots}  \hfill & {\ldots}  \hfill & {\ldots}  \hfill
\\
 {L_{1\,n}}  \hfill & {L_{2\,n}}  \hfill & {\ldots}  \hfill & {L_{n\,n}}  \hfill
\\
\end{array}} } \right),
\]
 where $L_{ij} $ $(\forall i,j=\overline{1,n})$ is a cofactor
  in $\lambda {\rm {\bf
I}} + {\rm {\bf A}}^{ *} {\rm {\bf A}}$.   By Lemma \ref{kyrc2},
${\rm {\bf A}}^{ +}  = {\mathop {\lim }\limits_{\lambda \to 0}}
\left( {\lambda {\rm {\bf I}} + {\rm {\bf A}}^{ * }{\rm {\bf A}}}
\right)^{ - 1}{\rm {\bf A}}^{ *} $, so that
\begin{equation}
\label{eq3} {\rm {\bf A}}^{ +}  ={\mathop {\lim} \limits_{\lambda
\to 0}}
\begin{pmatrix}
 \frac{\det \left( {\lambda {\rm {\bf I}} + {\rm {\bf A}}^{ *} {\rm {\bf A}}}
\right)_{.1} \left( {{\rm {\bf a}}_{.\,1}^{ *} }  \right)}{{\det
\left( {\lambda {\rm {\bf I}} + {\rm {\bf A}}^{ *} {\rm {\bf A}}}
\right)}}
 & \ldots & \frac{\det \left( {\lambda {\rm {\bf I}} + {\rm {\bf
A}}^{ *} {\rm {\bf A}}} \right)_{.\,1} \left( {{\rm {\bf
a}}_{.\,m}^{ *} }  \right)}{{\det \left( {\lambda {\rm {\bf I}} +
{\rm {\bf
A}}^{ *} {\rm {\bf A}}} \right)}}\\
  \ldots & \ldots & \ldots \\
\frac{\det \left( {\lambda {\rm {\bf I}} + {\rm {\bf A}}^{ *} {\rm
{\bf A}}} \right)_{.\,n} \left( {{\rm {\bf a}}_{.\,1}^{ *} }
\right)}{{\det \left( {\lambda {\rm {\bf I}} + {\rm {\bf A}}^{ *}
{\rm {\bf A}}} \right)}} &
 \ldots &
 \frac{\det \left( {\lambda {\rm {\bf I}} + {\rm {\bf A}}^{ *} {\rm {\bf
A}}} \right)_{.\,n} \left( {{\rm {\bf a}}_{.\,m}^{ *} }
\right)}{{\det \left( {\lambda {\rm {\bf I}} + {\rm {\bf A}}^{ *}
{\rm {\bf A}}} \right)}}
\end{pmatrix}.
\end{equation}
 From Theorem \ref{kyrc4} we get
\[ \det \left( {\lambda{\rm {\bf I}} + {\rm {\bf A}}^{ *}
{\rm {\bf A}}} \right) = \lambda ^{n} + d_{1} \lambda ^{n - 1} +
d_{2} \lambda ^{n - 2} + \ldots + d_{n},
\]
 where $d_{r} $ $(\forall r=\overline{1,n-1})$ is a
sum of  principal minors  of ${\rm {\bf A}}^{ *} {\rm {\bf A}}$ of
order $r$  and $d_{n}=\det  {\rm {\bf A}}^{ *} {\rm {\bf A}}$.
Since $\rank{\rm {\bf A}}^{ *} {\rm {\bf A}} = \rank{\rm {\bf A}}
= r$, then $d_{n} = d_{n - 1} = \ldots = d_{r + 1} = 0$ and
\begin{equation}
\label{eq4}\det \left( {\lambda {\rm {\bf I}} + {\rm {\bf A}}^{ *}
{\rm {\bf A}}} \right) = \lambda ^{n} + d_{1} \lambda ^{n - 1} +
d_{2} \lambda ^{n - 2} + \ldots + d_{r} \lambda ^{n - r}.
 \end{equation}
In the same way,  we have for  arbitrary $ 1\leq i \leq n$ and
$1\leq j\leq m $ from Theorem \ref{kyrc4}
\[ \det \left( {\lambda {\rm {\bf I}} + {\rm {\bf
A}}^{ *} {\rm {\bf A}}} \right)_{.\,i} \left( {{\rm {\bf
a}}_{.j}^{ *} } \right) = l_{1}^{\left( {ij} \right)} \lambda ^{n
- 1} + l_{2}^{\left( {ij} \right)} \lambda ^{n - 2} + \ldots +
l_{n}^{\left( {ij} \right)}, \]
 where for an arbitrary $1\leq k \leq n - 1$,\,
 $l_{k}^{\left( {ij}
\right)} = {\sum\limits_{\beta \in J_{k,\,n} {\left\{ {i}
\right\}}} {{\left| {\left( {({\rm {\bf A}}^{ *} {\rm {\bf
A}})_{.\,i} ({\rm {\bf a}}_{.j}^{ *}  )} \right)_{\beta} ^{\beta}
}  \right|}}}$, and \,\,$l_{n}^{\left( {i\,j} \right)} = \det
\left( {{\rm {\bf A}}^{ *} {\rm {\bf A}}} \right)_{.\,i} \left(
{{\rm {\bf a}}_{.\,j}^{ *} } \right)$. By Lemma \ref{kyrc5},
$\rank\left( {{\rm {\bf A}}^{ *} {\rm {\bf A}}} \right)_{.\,i}
\left( {{\rm {\bf a}}_{.\,j}^{ *} } \right) \le r$ so that if $
k>r$, then ${{\left| {\left( {({\rm {\bf A}}^{ *} {\rm {\bf
A}})_{\,.\,i} ({\rm {\bf a}}_{.j}^{ *}  )} \right)_{\beta}
^{\beta} } \right|}}= 0$, $(\forall \beta \in J_{k,\,n} {\left\{
{i} \right\}}, \forall i = \overline {1,n}, \forall j = \overline
{1,m})$. Therefore if $r + 1 \le k < n$, then $l_{k}^{\left( {ij}
\right)} = {\sum\limits_{\beta \in J_{k,\,n} {\left\{ {i}
\right\}}} {{\left| {\left( {({\rm {\bf A}}^{ *} {\rm {\bf
A}})_{\,.\,i} ({\rm {\bf a}}_{.j}^{ *}  )} \right)_{\beta}
^{\beta} } \right|}}}= 0$ and $l_{n}^{\left( {i\,j} \right)} =
\det \left( {{\rm {\bf A}}^{ *} {\rm {\bf A}}} \right)_{.\,i}
\left( {{\rm {\bf a}}_{.\,j}^{ *} } \right) = 0$, $\left( {\forall
i = \overline {1,n},\, \forall j = \overline {1,m}} \right)$.
Finally we obtain
\begin{equation}\label{eq5}
  \det \left( {\lambda {\rm {\bf I}} + {\rm {\bf A}}^{ *} {\rm {\bf
A}}} \right)_{.\,i} \left( {{\rm {\bf a}}_{.\,j}^{ *} } \right) =
l_{1}^{\left( {i\,j} \right)} \lambda ^{n - 1} + l_{2}^{\left(
{i\,j} \right)} \lambda ^{n - 2} + \ldots + l_{r}^{\left( {ij}
\right)} \lambda ^{n - r}.
\end{equation}

 By replacing the denominators and the numerators of the fractions in entries
 of  matrix (\ref{eq3}) with  the
 expressions (\ref{eq4}) and (\ref{eq5}) respectively, we get
\[
   {\rm {\bf A}}^{ +}  = {\mathop {\lim} \limits_{\lambda \to 0}}
\begin{pmatrix}
 {{\frac{{l_{1}^{\left( {11} \right)} \lambda ^{n - 1} + \ldots +
l_{r}^{\left( {11} \right)} \lambda^{n - r}}}{{\lambda ^{n} +
d_{1} \lambda ^{n - 1} + \ldots + d_{r} \lambda ^{n - r}}}}}
 & \ldots & {{\frac{{l_{1}^{\left( {1m} \right)} \lambda ^{n
- 1} + \ldots + l_{r}^{\left( {1m} \right)} \lambda ^{n -
r}}}{{\lambda ^{n} + d_{1} \lambda
^{n - 1} + \ldots + d_{r} \lambda ^{n - r}}}}} \\
  \ldots & \ldots & \ldots \\
  {{\frac{{l_{1}^{\left( {n1} \right)} \lambda ^{n - 1} + \ldots +
l_{r}^{\left( {n1} \right)} \lambda ^{n - r}}}{{\lambda ^{n} +
d_{1} \lambda ^{n - 1} + \ldots + d_{r} \lambda ^{n - r}}}}} &
\ldots & {{\frac{{l_{1}^{\left( {nm} \right)} \lambda ^{n - 1} +
\ldots + l_{r}^{\left( {nm} \right)} \lambda ^{n - r}}}{{\lambda
^{n} + d_{1} \lambda ^{n - 1} + \ldots + d_{r} \lambda ^{n -
r}}}}}
\end{pmatrix}=\]
  \[ = \left( {{\begin{array}{*{20}c}
 {{\frac{{l_{r}^{\left( {11} \right)}} }{{d_{r}} }}} \hfill & {\ldots}
\hfill & {{\frac{{l_{r}^{\left( {1m} \right)}} }{{d_{r}} }}} \hfill \\
 {\ldots}  \hfill & {\ldots}  \hfill & {\ldots}  \hfill \\
 {{\frac{{l_{r}^{\left( {n1} \right)}} }{{d_{r}} }}} \hfill & {\ldots}
\hfill & {{\frac{{l_{r}^{\left( {nm} \right)}} }{{d_{r}} }}} \hfill \\
\end{array}} } \right).
\]
From here  the representation (\ref{eq1}) of ${\rm {\bf A}}^{ +} $
follows by denoting $l_{r}^{\left( {ij} \right)} = l_{ij} $.

 We obtain the representation (\ref{eq2}) in the same way. $\blacksquare$
 \newtheorem{Remark}{Remark}[section]
\begin{Remark}\label{rem1} If $\rank{\rm {\bf A}} =
n$, then   from   Lemma \ref{kyrc3} we get ${\rm {\bf A}}^{ +}  =
\left( {{\rm {\bf A}}^{ *} {\rm {\bf A}}} \right)^{ - 1}{\rm {\bf
A}}^{ *} $. Representing $({\rm {\bf A}}^{ *} {\rm {\bf A}})^{-1}$
by the classical adjoint matrix, we  have
 \begin{equation}
\label{eq6} {\rm {\bf A}}^{ +}  = {\frac{{1}}{{\det ({\rm {\bf
A}}^{ *} {\rm {\bf A}})}}}
\begin{pmatrix}
   {\det ({\rm {\bf A}}^{ *} {\rm {\bf A}})_{.1} \left( {{\rm {\bf a}}_{.1}^{ *} }  \right)}&
   \ldots & {\det ({\rm
{\bf A}}^{ *} {\rm {\bf A}})_{.1} \left( {{\rm {\bf
a}}_{.\;m}^{ *} }  \right)}\\
  \ldots & \ldots & \ldots \\
  {\det ({\rm {\bf A}}^{ *} {\rm {\bf A}})_{.\,n} \left( {{\rm {\bf a}}_{.\,1}^{ *} }  \right)}
   & \ldots
   & {\det ({\rm {\bf A}}^{ *} {\rm {\bf A}})_{.\,n} \left( {{\rm {\bf a}}_{.\,m}^{ *} }  \right)}
\end{pmatrix}.\end{equation}
 If $n < m$, then  (\ref{eq2}) is
 valid.
 \end{Remark}
 \begin{Remark}\label{rem2}
 As above, if $\rank{\rm {\bf A}} = m$, then
\begin{equation}
\label{eq7} {\rm {\bf A}}^{ +}  = {\frac{{1}}{{\det ({\rm {\bf
A}}{\rm {\bf A}}^{ *} )}}}\begin{pmatrix}
  {\det ({\rm {\bf A}}{\rm {\bf A}}^{ *} )_{1\,.} \left( {{\rm {\bf a}}_{1\,.}^{ *} }  \right)}
   & \ldots &
   {\det ({\rm {\bf A}}{\rm {\bf A}}^{ *} )_{m\,.} \left( {{\rm {\bf a}}_{1\,.}^{ *} }  \right)} \\
 \ldots & \ldots & \ldots \\
 {\det ({\rm {\bf A}}{\rm {\bf A}}^{ *} )_{1\,.} \left( {{\rm {\bf a}}_{n\,.}^{ *} }  \right)}
  & \ldots &
{\det ({\rm {\bf A}}{\rm {\bf A}}^{ *} )_{m\,.} \left( {{\rm {\bf
a}}_{n\,.}^{ *} }  \right)}
\end{pmatrix}.
\end{equation}
If  $n > m$, then  (\ref{eq1}) is valid as well.
 \end{Remark}
 \begin{Remark}
The representation (\ref{eq1}) can be obtained from Theorem
\ref{kyrc1} by using the Binet-Cauchy formula. We use another
method, which is the same for the determinantal representations of
the Moore-Penrose inverse by  (\ref{eq1}) and (\ref{eq2}), and of
the Drazin inverse by (\ref{eq11}).
 \end{Remark}
 \begin{Remark} To obtain an entry of ${\rm {\bf A}}^{ +}$ by
 Theorem \ref{kyrc1} one calculates
 $(C^{r}_{n}C^{r}_{m}+C^{r-1}_{n-1}C^{r-1}_{m-1})$ determinants of
 order $r$. Whereas by  (\ref{eq1}) we calculate as much
 as
 $(C^{r}_{n}+C^{r-1}_{n-1})$ determinants of
 order $r$ or  we calculate the total of
 $(C^{r}_{m}+C^{r-1}_{m-1})$ determinants by (\ref{eq2}). Therefore the
 calculation of entries of ${\rm {\bf A}}^{ +} $ by Theorem \ref{kyrc6} is easier
  than by Theorem \ref{kyrc1}.
  \end{Remark}
\newtheorem{Corollary}{Corollary}[section]
 \begin{Corollary}\label{kyrc7}
If ${\rm {\bf A}}\in {\mathbb C}^{m\times n}_{r}$ and $r < \min
{\left\{ {m,n} \right\}}$  or $r = m < n $, then the projection
matrix ${\rm {\bf P}} = {\rm {\bf A}}^{ +} {\rm {\bf A}}$ can be
represented as
\[
{\rm {\bf P}} = \left({\frac{{p_{ij}}}{{d_{r} \left( {{\rm {\bf
A}}^{ *} {\rm {\bf A}}} \right)}}}\right)_{n\times n},
\]
 where ${\rm {\bf d}}_{.\,j} $ denotes the
$j$th column of $({\rm {\bf A}}^{ *} {\rm {\bf A}})$ and, for
arbitrary $1\leq i,j \leq n $,
   $p_{ij} ={\sum\limits_{\beta \in
J_{r,n} {\left\{ {i} \right\}}} {{\left| {\left( {({\rm {\bf A}}^{
*} {\rm {\bf A}})_{\,.\,i} ({\rm {\bf d}}_{.j}  )} \right)_{\beta}
^{\beta} }  \right|}}}$.
\end{Corollary}
{\textit{Proof}}. Representing the Moore - Penrose inverse ${\rm
{\bf A}}^{+}$ by (\ref{eq1}), we obtain
\[
{\rm {\bf P}}={\frac{{1}}{{d_{r} \left( {{\rm {\bf A}}^{ *} {\rm
{\bf A}}} \right)}}}\begin{pmatrix}
  l_{11} & l_{12} & \ldots & l_{1m} \\
  l_{21} & l_{22} & \ldots & l_{2m} \\
  \ldots & \ldots & \ldots & \ldots \\
  l_{n1} & l_{n2} & \ldots & l_{nm}
\end{pmatrix}\begin{pmatrix}
  a_{11} & a_{12} & \ldots & a_{1n} \\
  a_{21} & a_{22} & \ldots & a_{2n} \\
  \ldots & \ldots & \ldots & \ldots \\
  a_{m\,1} & a_{m\,2} & \ldots & a_{m\,n}
\end{pmatrix}.
\]
Therefore, for arbitrary $1\leq i,j\leq n$ we get
 \[\begin{array}{c}
  p_{i\,j}=\sum\limits_{k}
 {{\sum\limits_{\beta \in J_{r,\,
n} {\left\{ {i} \right\}}} {{\left| {\left( {({\rm {\bf A}}^{ *}
{\rm {\bf A}})_{.\,i} ({\rm {\bf a}}_{.\,k}^{ *}  )}
\right)_{\beta} ^{\beta} }  \right|}}}} \cdot a_{k\,j}=\\
  =\sum\limits_{\beta \in J_{r,\,
n} {\left\{ {i} \right\}}} \sum\limits_{k}
 {{ {{\left|
 {\left( {({\rm {\bf A}}^{ *}
{\rm {\bf A}})_{.\,i} ({\rm {\bf a}}_{.\,k}^{ *}\cdot a_{k\,j}  )}
\right)_{\beta} ^{\beta} } \right|}}}}={\sum\limits_{\beta \in
J_{r,\,n} {\left\{ {i} \right\}}} {{\left| {\left( {({\rm {\bf
A}}^{ *} {\rm {\bf A}})_{.\,i} ({\rm {\bf d}}_{.j}^{ *}  )}
\right)_{\beta} ^{\beta} } \right|}}}.\blacksquare
\end{array}
\]

Using the representation (\ref{eq2}) of the Moore - Penrose
inverse the following corollary can be proved in the same way.
 \begin{Corollary} \label{kyrc8} If ${\rm
{\bf A}}\in {\mathbb C}^{m\times n}_{r}$, where $r < \min {\left\{
{m,n} \right\}}$ or $r = n < m$,
 then a projection matrix ${\rm {\bf Q}} = {\rm
{\bf A}}{\rm {\bf A}}^{ +} $  can be represented as
\[
{\rm {\bf Q}} =  \left({\frac{{q_{ij}}}{{d_{r} \left( {{\rm {\bf
A}}{\rm {\bf A}}^{ *} } \right)}}}\right)_{m\times m},
\]
where ${\rm {\bf g}}_{i.} $ denotes the $i$th row of $({\rm {\bf
A}}{\rm {\bf A}}^{ *} )$ and,  for arbitrary $1\leq i,j \leq m $,
$q_{i\,j} ={\sum\limits_{\alpha \in I_{r,m} {\left\{ {j}
\right\}}} {{\left| {\left( {({\rm {\bf A}} {\rm {\bf A}}^{
*})_{j.}\, ({\rm {\bf g}}_{i.}\,  )} \right)_{\alpha} ^{\alpha} }
\right|}}} $.
\end{Corollary}

 \begin{Remark} By definition of the classical adjoint $Adj({\rm {\bf A}})$
  for an arbitrary invertible matrix ${\rm {\bf A}}\in {\mathbb C}^{n\times
  n}$ one may put $Adj({\rm {\bf A}})\cdot{\rm {\bf A}}=\det{\rm {\bf A}}\cdot{\rm {\bf I}}_{n}$.
  If ${\rm {\bf A}}\in {\mathbb C}^{m\times n}$ and $\rank{\rm {\bf
A}} = n$, then by  Lemma \ref{kyrc3}, ${\rm {\bf A}}^{ +} {\rm
{\bf A}} = {\rm {\bf I}}_{n}$.
  Representing the matrix ${\rm {\bf A}}^{ +} $
 by  (\ref{eq6}) as ${\rm
{\bf A}}^{ +}  = {\frac{{{\rm {\bf L}}}}{{\det \left( {{\rm {\bf
A}}^{ *} {\rm {\bf A}}} \right)}}}$, we obtain
 ${\rm {\bf L}}{\rm {\bf A}} = \det \left(
{{\rm {\bf A}}^{ *} {\rm {\bf A}}} \right) \cdot {\rm {\bf I}}_{n}
$. This means that the matrix ${\rm {\bf L}}$
 is a
 left  analogue of $Adj({\rm {\bf A}})$ , where ${\rm {\bf A}}\in {\mathbb C}^{m\times n}_{n}$.
  If $\rank{\rm {\bf A}} = m$, then by  Lemma \ref{kyrc3}, ${\rm
{\bf A}}{\rm {\bf A}}^{ +} = {\rm {\bf I}}_{m}$. Representing the
matrix ${\rm {\bf A}}^{ +} $
 by  (\ref{eq7}) as ${\rm {\bf A}}^{ +}  = {\frac{{{\rm
{\bf R}}}}{{\det \left( {{\rm {\bf A}}{\rm {\bf A}}^{ *} }
\right)}}}$, we obtain ${\rm {\bf A}}{\rm {\bf R}} = {\rm {\bf
I}}_{m} \cdot \det \left( {{\rm {\bf A}}{\rm {\bf A}}^{ *} }
\right)$. This means that the matrix ${\rm {\bf R}}$ is a
 right analogue of $Adj({\rm {\bf A}})$, where ${\rm {\bf A}}\in {\mathbb C}^{m\times n}_{m}$.
\end{Remark}
\begin{Remark}\label{kyrc9}
  If
${\rm} {\rm {\bf A}}\in {\mathbb C}^{m\times n}_{r} $ and
$r<\min\{m, n\}$, then by (\ref{eq1}) we have ${\rm {\bf A}}^{ +}
= {\frac{{{\rm {\bf L}}}}{{d_{r} \left( {{\rm {\bf A}}^{ *} {\rm
{\bf A}}} \right)}}}$,
 where ${\rm {\bf L}} = \left( {l_{ij}}  \right)\in {\mathbb C}^{n\times m}$.
  From Corollary \ref{kyrc7}  we get ${\rm
{\bf L}}{\rm {\bf A}} = d_{r} \left( {\rm {\bf A}}^{ *} {\rm {\bf
A}} \right) \cdot {\rm {\bf P}}$.  The matrix ${\rm {\bf P}}$ is
idempotent. All eigenvalues of an idempotent matrix chose from 1
or 0 only. Thus, there exists an unitary matrix ${\rm {\bf U}}$
such that ${\rm {\bf L}}{\rm {\bf A}} = d_{r} \left( {\rm {\bf
A}}^{ *} {\rm {\bf A}} \right) {\rm {\bf U}} {{\rm {\bf
diag}}\left( {{\rm 1},\ldots,{\rm 1},{\rm 0},\ldots ,{\rm 0}}
\right)}  {\rm {\bf U}}^{ *} $, where $
 {{\rm {\bf diag}}\left( {{\rm 1},\ldots,{\rm 1},{\rm
0},\ldots ,{\rm 0}} \right)} \in {\mathbb C}^{n\times n}$ is a
diagonal matrix.  Therefore,  the matrix ${\rm {\bf L}}$
 can be considered  as a left
analogue of $Adj({\rm {\bf A}})$, where ${\rm {\bf A}}\in {\mathbb
C}^{m\times n}_{r}$.

 In the same
way, if ${\rm} {\rm {\bf A}}\in {\mathbb C}^{m\times n}_{r} $ and
$r<\min\{m, n\}$, then by   (\ref{eq2}) we have ${\rm {\bf A}}^{
+}  = {\frac{{{\rm {\bf R}}}}{{d_{r} \left( {{\rm {\bf A}}{\rm
{\bf A}}^{ *} } \right)}}}$, where ${\rm {\bf R}} = \left(
{r_{ij}} \right)\in {\mathbb C}^{n\times m}$. From Corollary
\ref{kyrc8}  we get  ${\rm {\bf A}}{\rm {\bf R}} = d_{r} \left(
{\rm {\bf A}}{\rm {\bf A}}^{ *} \right) \cdot {\rm {\bf Q}}$. The
matrix ${\rm {\bf Q}}$ is idempotent. There exists an unitary
matrix ${\rm {\bf V}}$ such that ${\rm {\bf A}}{\rm {\bf R}} =
d_{r} \left( {\rm {\bf A}} {\rm {\bf A}}^{ *} \right) {\rm {\bf
V}} {{\rm {\bf diag}}\left( {{\rm 1},\ldots ,{\rm 1},{\rm
0},\ldots,{\rm 0}} \right)}{\rm {\bf V}}^{ *} $, where $ {{\rm
{\bf diag}}\left( {{\rm 1},\ldots ,{\rm 1},{\rm 0},\ldots,{\rm 0}}
\right)} \in {\mathbb C}^{m\times m}$. Therefore,  the matrix
${\rm {\bf R}}$ can be considered  as a right analogue of
$Adj({\rm {\bf A}})$ in this case.
\end{Remark}

\section{An analogue of the classical adjoint matrix for the Drazin
inverse}
\begin{definition} \cite{ca,dr}
Let  ${\rm {\bf A}}\in {\mathbb C}^{n\times n} $ with $ Ind{\kern
1pt} {\rm {\bf A}}=k$, where  a nonnegative integer $Ind{\kern
1pt} {\rm {\bf A}}:= {\mathop {\min} \limits_{k \in N \cup
{\left\{ {0} \right\}}} }{\kern 1pt} {\left\{ {\rank{\rm {\bf
A}}^{k + 1} = \rank{\rm {\bf A}}^{k}} \right\}}$.
 Then the matrix ${\rm {\bf X}}$
satisfying \begin{equation} \label{eq8} {\rm {\bf A}}^{k+1}{\rm
{\bf X}}={\rm {\bf
         A}}^{k};\,\,
{\rm {\bf X}}{\rm {\bf A}}{\rm {\bf X}}={\rm {\bf X}};\,\, {\rm
{\bf A}}{\rm {\bf X}}={\rm {\bf X}}{\rm {\bf A}}
\end{equation}
 is called the Drazin
inverse of ${\rm {\bf A}}$ and is denoted by ${\rm {\bf X}}={\rm
{\bf A}}^{D}$. In particular, if $Ind{\kern 1pt} {\rm {\bf A}}=1$,
then the matrix ${\rm {\bf X}}$ in (\ref{eq8}) is called the group
inverse and is denoted by ${\rm {\bf X}}={\rm {\bf A}}^{{\rm \#}
}$.
\end{definition}
\begin{Remark}If $Ind{\kern 1pt}
{\rm {\bf A}}=0$, then ${\rm {\bf A}}$ is nonsingular, and ${\rm
{\bf A}}^{D}\equiv {\rm {\bf A}}^{-1}$.
\end{Remark}
The Drazin inverse can be represented explicitly by the Jordan
canonical form as follows.
\begin{theorem}\cite{ca} If ${\rm {\bf A}} \in {\mathbb C}^{n\times n}$ with
$Ind{\kern 1pt} {\rm {\bf A}} = k $ and
\begin{equation}\label{eq9}
{\rm {\bf A}} = {\rm {\bf P}}\begin{pmatrix}
  {\rm {\bf C}} & {\rm {\bf 0}} \\
  {\rm {\bf 0}}& {\rm {\bf N}}
\end{pmatrix} {\rm {\bf P}}^{ - 1},
\end{equation}
where ${\rm {\bf C}}$ is nonsingular, $\rank{\rm {\bf C}} =
\rank{\rm {\bf A}}^{k}$, and ${\rm {\bf N}}$ is nilpotent of order
$k$, then
\[
{\rm {\bf A}}^{D} = {\rm {\bf P}}\begin{pmatrix}
  {\rm {\bf C}}^{ - 1} & {\rm {\bf 0}} \\
  {\rm {\bf 0}}& {\rm {\bf 0}}
\end{pmatrix} {\rm {\bf P}}^{ - 1}.
\]
\end{theorem}

We use the following theorem about the limit representation of the
Drazin inverse.
\begin{theorem} \cite{ca}\label{kyrc10} If ${\rm {\bf A}}
 \in {\mathbb C}^{n\times n}$, then
\[
{\rm {\bf A}}^{D} = {\mathop {\lim} \limits_{\lambda \to 0}}
\left( {\lambda {\rm {\bf I}}_n + {\rm {\bf A}}^{k + 1}} \right)^{
- 1}{\rm {\bf A}}^{k},
\]
where $k = Ind{\kern 1pt} {\rm {\bf A}}$ and $\lambda \in {\mathbb
R} _{ +}  $.
\end{theorem}
Denote by ${\rm {\bf a}}_{.j}^{(k)} $ and ${\rm {\bf
a}}_{i.}^{(k)} $ the $j$th column  and the $i$th row of  ${\rm
{\bf A}}^{k} $ respectively.
\begin{lemma} \label{kyrc11} If ${\rm {\bf A}} \in {\mathbb C}^{n\times n}$
with $Ind{\kern 1pt} {\rm {\bf A}} = k $, then
\begin{equation}
\label{eq10} \rank{\rm {\bf A}}_{.{\kern 1pt} i}^{k + 1} \left(
{{\rm {\bf a}}_{.j}^{\left( {k} \right)}}  \right) \le \rank{\rm
{\bf A}}^{k + 1}, \quad \forall i,j = \overline {1,n}.
\end{equation}
\end{lemma}
 {\textit{Proof}}. The proof of this lemma is similar to that of
 Lemma \ref{kyrc5}.

 In the following theorem we introduce a determinantal
 representation of the Drazin inverse.
\begin{theorem}\label{kyrc12}
 If $Ind{\kern 1pt} {\rm {\bf A}} = k $ and $\rank{\rm {\bf
A}}^{k + 1} = \rank{\rm {\bf A}}^{k}=r \le n$ for an arbitrary
matrix ${\rm {\bf A}}\in {\mathbb C}^{n\times n} $, then
\begin{equation}
\label{eq11}
\begin{array}{l}
  {\rm {\bf A}}^{D} = \left({\frac{{d_{ij}}}{{d_{r} \left( {{\rm
{\bf A}}^{k + 1}} \right)}}}\right)_{n\times n},
  \\
 \mbox{where}\\
  d_{r} \left( {{\rm {\bf A}}^{k + 1}} \right) =
{\sum\limits_{\beta \in J_{r,n}} {{\left| {\left( {{\rm {\bf
A}}^{k + 1}} \right)_{\beta} ^{\beta }}  \right|}}},\\
d_{ij} = {\sum\limits_{\beta \in J_{r,n} {\left\{ {i} \right\}}}
{{\left| {\left( {{\rm {\bf A}}_{.\,i}^{k + 1} \left( {{\rm {\bf
a}}_{.j}^{\left( {k} \right)}}  \right)} \right)_{\beta} ^{\beta}
} \right|}}}, \,\,\left( {\forall i,j = \overline {1,n}} \right).
\end{array}
\end{equation}
\end{theorem}
{\textit{Proof}}. The proof of this theorem is analogous to that
of Theorem \ref{kyrc6} by using Theorem \ref{kyrc4}, Lemma
\ref{kyrc11}, and Theorem \ref{kyrc10}.

In the following corollaries we introduce  determinantal
 representations of the group inverse ${\rm {\bf A}}^{{\rm \#} }$ and
  the matrix ${\rm {\bf A}}^{D}{\rm {\bf A}}$ respectively.
\begin{Corollary}
 If $Ind{\kern 1pt} {\rm {\bf A}} = 1 $ and $\rank{\rm {\bf
A}}^{2} = \rank{\rm {\bf A}}=r \le n$ for an arbitrary matrix
${\rm {\bf A}}\in {\mathbb C}^{n\times n} $, then
\[
  {\rm {\bf A}}^{\rm \#} = \left({\frac{{g_{ij}}}{{d_{r} \left( {{\rm
{\bf A}}^{2}} \right)}}}\right)_{n\times n},
  \]
where
  $d_{r} \left( {{\rm {\bf A}}^{2}} \right) =
{\sum\limits_{\beta \in J_{r,n}} {{\left| {\left( {{\rm {\bf
A}}^{2}} \right)_{\beta} ^{\beta }}  \right|}}}$, $ g_{ij} =
{\sum\limits_{\beta \in J_{r,n} {\left\{ {i} \right\}}} {{\left|
{\left( {{\rm {\bf A}}_{.\,i}^{2} \left( {{\rm {\bf a}}_{.j}}
\right)} \right)_{\beta} ^{\beta} } \right|}}}$, $\left( {\forall
i,j = \overline {1,n}} \right). $
\end{Corollary}
{\textit{Proof}}. The proof follows  from Theorem \ref{kyrc12} in
view of $k=1$.
\begin{Corollary}
If $Ind{\kern 1pt} {\rm {\bf A}} = k $ and $\rank{\rm {\bf A}}^{k
+ 1} = \rank{\rm {\bf A}}^{k}=r \le n$ for an arbitrary matrix
${\rm {\bf A}}\in {\mathbb C}^{n\times n} $, then
\[
  {\rm {\bf A}}^{D}{\rm {\bf A}} = \left({\frac{{v_{ij}}}{{d_{r} \left( {{\rm
{\bf A}}^{k+1}} \right)}}}\right)_{n\times n},
  \]
where
 $ v_{ij} =
{\sum\limits_{\beta \in J_{r,n} {\left\{ {i} \right\}}} {{\left|
{\left( {{\rm {\bf A}}_{.\,i}^{k+1} \left( {{\rm {\bf
a}}_{.j}}^{(k+1)} \right)} \right)_{\beta} ^{\beta} } \right|}}}$,
$\left( {\forall i,j = \overline {1,n}} \right). $
\end{Corollary}
{\textit{Proof}}. Representing the Drazin inverse ${\rm {\bf
A}}^{D}$ by (\ref{eq11}) we obtain
\[
{\rm {\bf A}}^{D}{\rm {\bf A}}= \left({\frac{{d_{ij}}}{{d_{r}
\left( {{\rm {\bf A}}^{k+1}} \right)}}}\right)_{n\times n}\cdot
\left( a_{ij}\right)_{n\times n}= \left({\frac{{v_{ij}}}{{d_{r}
\left( {{\rm {\bf A}}^{k+1}} \right)}}}\right)_{n\times n}.
\]
Here for arbitrary $1\leq i,j\leq n$ we have
 \[\begin{array}{c}
  v_{i\,j}=\sum\limits_{s}
 {{\sum\limits_{\beta \in J_{r,\,
n} {\left\{ {i} \right\}}} {{\left| {\left( {({\rm {\bf A}}^{
k+1})_{.\,i} ({\rm {\bf a}}_{.\,s}^{ (k)}  )}
\right)_{\beta} ^{\beta} }  \right|}}}} a_{s\,j}=\\
  =\sum\limits_{\beta \in J_{r,\,
n} {\left\{ {i} \right\}}} \sum\limits_{s}
 {{ {{\left|
 {\left( {({\rm {\bf A}}^{k+1}
)_{.\,i} ({\rm {\bf a}}_{.\,s}^{(k)} \cdot a_{s\,j} )}
\right)_{\beta} ^{\beta} } \right|}}}}={\sum\limits_{\beta \in
J_{r,\,n} {\left\{ {i} \right\}}} {{\left| {\left( {({\rm {\bf
A}}^{ k+1} )_{.\,i} ({\rm {\bf a}}_{.j}^{(k+1)}  )}
\right)_{\beta} ^{\beta} } \right|}}}. \blacksquare
\end{array}
\]
\begin{Remark} The matrix $({\rm {\bf A}}^{D}{\rm {\bf A}})$ is idempotent. Similarly
 to the case of Remark \ref{kyrc9}, the matrix ${\rm {\bf
D}}=\left(d_{ij}\right)_{n\times n}$ can be considered as an
analogue of the classical adjoint matrix, where ${\rm {\bf
D}}={d_{r} \left( {{\rm {\bf A}}^{k+1}} \right)}{\rm {\bf
A}}^{D}$.
\end{Remark}

\section{ Cramer rules for generalized inverse solutions}
\begin{definition} Suppose in a complex system of  linear
equations:
\begin{equation}
\label{eq12} {\rm {\bf A}} \cdot {\rm {\bf x}} = {\rm {\bf y}}
\end{equation}
the coefficient matrix ${\rm {\bf A}} \in {\mathbb C}^{m\times
n}_{r}$ and a column of constants ${\rm {\bf y}} = \left( {y_{1}
,\ldots ,y_{m} } \right)^{T}\in {\mathbb C}^{m}$. The least
squares solution of the system (\ref{eq12})
 is  the vector ${\rm {\bf x}}^{0}\in\mathbb{C}^{n} $
  satisfying
 \[{\left\| {{\rm
{\bf x}}^{0}} \right\|} = {\mathop {\min} \limits_{{\rm {\bf
\tilde {x}}} \in\mathbb{C}^{n}} } {\left\{ {{\left\| {{\rm {\bf
\tilde {x}}}} \right\|}\,\vert\, {\left\| {{\rm {\bf A}} \cdot
{\rm {\bf \tilde {x}}} - {\rm {\bf y}}} \right\|} = {\mathop
{\min} \limits_{{\rm {\bf x}} \in \mathbb{C}^{n}} }{\left\| {{\rm
{\bf A}} \cdot {\rm {\bf x}} - {\rm {\bf y}}} \right\|}}
\right\}},\] where $\mathbb{C}^{n}$ is an $n$-dimension complex
vector space.
\end{definition}
 \begin{theorem}
\cite{ho}  The vector ${\rm {\bf x}} = {\rm {\bf A}}^{ +} {\rm
{\bf y}}$
 is the least squares
solution  of the system (\ref{eq12}).
\end{theorem}
\begin{theorem}\label{kyrc13}
 The following
statements are true for the system of  linear equations
(\ref{eq12}).
\begin{itemize}
\item [ i)] If $\rank{\rm {\bf A}} = n$,
 then the components of the least squares
solution ${\rm {\bf x}}=\left( {{x^{0}_{1}} ,\ldots ,x^{0}_{n} }
\right)^{T}$ are obtained by the formula
\begin{equation}
\label{eq13} x_{j}^{0} = {\frac{{\det  ({\rm {\bf A}}^{ *} {\rm
{\bf A}})_{.\,j} \left( {{\rm {\bf f}}} \right)}}{{\det {\rm {\bf
A}}^{ *} {\rm {\bf A}}}}}, \quad \left( {\forall j = \overline
{1,n}} \right),
\end{equation}
 where  ${\rm {\bf f}} = {\rm {\bf A}}^{ *} {\rm {\bf y}}$.

 \item [ii)] If $\rank{\rm {\bf A}} = r \le m < n$, then
\begin{equation}
\label{eq14} x_{j}^{0} = {\frac{{\sum\limits_{\beta \in J_{r,n}
{\left\{ {j} \right\}}} {{\left| {\left( {({\rm {\bf A}}^{ *} {\rm
{\bf A}})_{.\,j} ({\rm {\bf f}}  )} \right)_{\beta} ^{\beta} }
\right|}}}}{{d_{r} \left( {{\rm {\bf A}}^{ *} {\rm {\bf A}}}
\right)}}}, \quad \left( {\forall j = \overline {1,n}} \right).
\end{equation}
\end{itemize}
\end{theorem}

{\textit{Proof.}} i) If $\rank{\rm {\bf A}} = n$, then  we can
represent ${\rm {\bf A}}^{ +} $ by (\ref{eq7}). By multiplying
${\rm {\bf A}}^{ +}$  into ${\rm {\bf y}}$ we get (\ref{eq13}).

ii) If $\rank{\rm {\bf A}} = k \le m < n$, then    ${\rm {\bf
A}}^{ +} $ can be  represented by  (\ref{eq1}). By multiplying
${\rm {\bf A}}^{ +}$  into ${\rm {\bf y}}$ the least squares
solution of the linear system (\ref{eq12}) is given by components
as in (\ref{eq14}). $\blacksquare$

Using (\ref{eq3}) and (\ref{eq8}), we can obtain another
representation of the Cramer rule for the least squares solution
of a
 linear system.
\begin{theorem}
 The following statements are
true for a system of  linear equations written in the form ${\rm
{\bf x}}\cdot {\rm {\bf A}} = {\rm {\bf y}}$.
\begin{itemize}
\item [ i)] If $\rank{\rm {\bf A}} = m$,
 then the components of the least squares
solution ${\rm {\bf x}}^{0}={\rm {\bf y}}{\rm {\bf A}}^{+}$ are
obtained by the formula
\[ x_{i}^{0} = {\frac{{\det  ({\rm {\bf A}}{\rm {\bf A}}^{ *} )_{i\,.}
\left( {{\rm {\bf g}}} \right)}}{{\det {\rm {\bf A}}{\rm {\bf
A}}^{ *} }}}, \quad \left( {\forall i = \overline {1,m}} \right),
\]
 where  ${\rm {\bf g}} = {\rm {\bf y}}{\rm {\bf A}}^{ *} $.

 \item [ii)] If $\rank{\rm {\bf A}} = r \le  n< m$, then
\[ x_{i}^{0} = {\frac{{\sum\limits_{\alpha \in I_{r,m} {\left\{ {i} \right\}}}
{{\left| {\left( {({\rm {\bf A}}{\rm {\bf A}}^{ *} )_{\,i\,.}
({\rm {\bf g}}  )} \right)_{\alpha} ^{\alpha} } \right|}}}}{{d_{r}
\left( {{\rm {\bf A}}{\rm {\bf A}}^{ *} } \right)}}}, \quad \left(
{\forall i = \overline {1,m}} \right).
\]
\end{itemize}
\end{theorem}
{\textit{Proof}}. The proof of this theorem is analogous to that
of Theorem \ref{kyrc13}.
\begin{Remark}
The obtained formulas of the Cramer rule for the least squares
solution differ from  similar formulas in \cite{ben2,ch,
ji,wa1,wa2,wei}. They give  a  closer approximation to the Cramer
rule for consistent nonsingular systems of linear equations.
\end{Remark}

 In some situations,
however, people pay more attention to the Drazin inverse solution
of  singular linear systems \cite{ch, wei}. Consider a general
system of linear equations (\ref{eq12}), where ${\rm {\bf A}} \in
{\mathbb C}^{n\times n}$ and ${\rm {\bf x}}$, ${\rm {\bf y}}$ are
 vectors in ${\mathbb C}^{n}$. $R({\rm {\bf A}})$ denotes the
range of ${\rm {\bf A}}$ and $N({\rm {\bf A}})$ denotes the null
space of ${\rm {\bf A}}$. The characteristic of the Drazin inverse
solution ${\rm {\bf A}}^{D}\rm {\bf y}$ is given in \cite{wei} by
the following theorem.

\begin{theorem}\label{kyrc14}
Let ${\rm {\bf A}} \in {\mathbb C}^{n\times n}$ with $Ind(A) = k$.
Then ${\rm {\bf A}}^{D}{\rm {\bf y}}$ is both the unique
solution in $R({\rm {\bf A}}^{k})$ of
\begin{equation}\label{eq15}
{\rm {\bf A}}^{k+1}\rm {\bf x} = {\rm {\bf A}}^{k}{\rm {\bf y}},
\end{equation}
and the unique minimal ${\rm {\bf P}}$-norm least squares solution
of (\ref{eq12}).
\end{theorem}

\begin{Remark}
 The ${\rm {\bf P}}$-norm is defined as $\|{\rm
{\bf x}}\|_{{\rm {\bf P}}} = \|{\rm {\bf P}}^{-1}{\rm {\bf x}}\|$
for ${\rm {\bf x}}\in {\mathbb C}^{n}$, where ${\rm {\bf P}}$ is a
nonsingular matrix that transforms ${\rm {\bf A}}$ into its Jordan
canonical form (\ref{eq9}).
\end{Remark}

\begin{Remark}
Since  (\ref{eq15}) is analogous to the normal system ${\rm {\bf
A}}^{\ast}{\rm {\bf A}}{\rm {\bf x}}={\rm {\bf A}}^{\ast}{\rm {\bf
y}}$, the system (\ref{eq15}) is called  the generalized normal
equations of (\ref{eq12}), (see  \cite{wei}).
\end{Remark}
We obtain the Cramer rule for the  ${\rm {\bf P}}$-norm least
squares solution of (\ref{eq12}) in the following theorem.
\begin{theorem}
Let ${\rm {\bf A}} \in {\mathbb C}^{n\times n}$ with $Ind({\rm
{\bf A}}) = k$. Then the unique minimal ${\rm {\bf P}}$-norm least
squares solution ${\rm {\bf \widehat{x}}}=(\widehat{x}_{1},
\ldots,\widehat{x}_{n})^{T}$ of the system (\ref{eq12}) is given
by
\begin{equation}
\label{eq16} \widehat{x}_{i} = {\frac{{{\sum\limits_{\beta \in
J_{r,\,n} {\left\{ {i} \right\}}} {{\left| {\left( {{\rm {\bf
A}}_{.{\kern 1pt} i}^{k + 1} \left( {{\rm {\bf g}}} \right)}
\right)_{\beta} ^{\beta} }  \right|}} }}}{{{\sum\limits_{\beta \in
J_{r,n}} {{\left| {\left( {{\rm {\bf A}}^{k + 1}} \right)_{\beta}
^{\beta} }  \right|}}} }}} \quad \forall i = \overline {1,n},
\end{equation}
where ${\rm {\bf g}} = {\rm {\bf A}}^{k}{\rm {\bf y}}.$
\end{theorem}
{\textit{Proof}}. Representing  the Drazin inverse by
(\ref{eq11}) and by virtue of Theorem \ref{kyrc14}, we have
\[
{\rm {\bf \widehat{x}}} =\begin{pmatrix}
  \widehat{x}_{1} \\
  \ldots \\
  \widehat{x}_{n}
\end{pmatrix}={\rm {\bf A}}^{D}{\rm {\bf y}}=\frac{1}{d_{r}
\left( {{\rm {\bf A}}^{k+1} } \right)}\begin{pmatrix}
  \sum\limits_{s = 1}^{n} {d_{1s} y_{s}}   \\
  \ldots \\
 \sum\limits_{s = 1}^{n} {d_{ns} y_{s}}
\end{pmatrix}.
\]
Therefore,
\[
\widehat{x}_{i} = {\frac{{1}}{{d_{r} \left( {{\rm {\bf A}}^{k +
1}} \right)}}}{\sum\limits_{s = 1\,\,}^{n} {{\sum\limits_{\beta
\in \, J_{r,n} {\left\{ {i} \right\}}} {{\left| {\left( {{\rm {\bf
A}}_{.{\kern 1pt} i}^{k + 1} \left( {{\rm {\bf a}}_{.\,s}^{\left(
{k} \right)}}  \right)} \right)_{\beta} ^{\beta} }  \right|}}} }}
\cdot y_{s} =
\]
\[
 = {\frac{{1}}{{d_{r} \left( {{\rm {\bf A}}^{k + 1}}
\right)}}}{\sum\limits_{\beta \in J_{r,\,n} {\left\{ {i}
\right\}}} {\,{\sum\limits_{\,s = 1\,\,}^{n} {{\left| {\left(
{{\rm {\bf A}}_{.\,i}^{k + 1} \left( {{\rm {\bf a}}_{.\,s}^{\left(
{k} \right)}}  \right)} \right)_{\beta} ^{\beta} }  \right|} \cdot
y_{s}} } }}=
\]
\[
 = {\frac{{1}}{{d_{r} \left( {{\rm {\bf A}}^{k + 1}}
\right)}}}{\sum\limits_{\beta \in J_{r,\,n} {\left\{ {i}
\right\}}} {\;{\sum\limits_{\,s = 1\,\,}^{n} {{\left| {\left(
{{\rm {\bf A}}_{.\,i}^{k + 1} \left( {{\rm {\bf a}}_{.\,s}^{\left(
{k} \right)} \cdot y_{s}}  \right)} \right)_{\beta} ^{\beta} }
\right|}}} }}.
\]
From this (\ref{eq16}) follows immediately. $\blacksquare$
\section{Examples}
1. Let us consider the system of linear equations.
\begin{equation} \label{eq17}
\left\{ \begin{array}{c}
  2x_1-5x_3+4x_4=1, \\
   7x_1-4x_2-9x_3+1.5x_4=2, \\
  3x_1-4x_2+7x_3-6.5x_4=3, \\
  x_1-4x_2+12x_3-10.5x_4=1.
\end{array}\right.\end{equation}
The coefficient matrix of the system is the matrix  ${\rm \bf
A}=\begin{pmatrix}
  2 & 0 & -5 & 4 \\
  7 & -4 & -9 & 1.5 \\
  3 & -4 & 7 & -6.5 \\
  1 & -4 & 12 & -10.5
\end{pmatrix}$. We calculate the rank of ${\rm {\bf A}}$ which is equal to 3, and we have
\[{\rm \bf A}^{\ast}=\begin{pmatrix}
  2 & 7 & 3 & 1 \\
  0 & -4 & -4 & -4 \\
  -5 & -9 & 7 & 12 \\
  4 & 1.5 & -6.5 & -10.5
\end{pmatrix}, {\rm \bf A}^{*}{\rm \bf A}=\begin{pmatrix}
  63 & -44 & -40 & -11.5 \\
-44 & 48 & -40 & 62 \\
  -40 & -40 & 299 & -205 \\
  -11.5 & 62 & -205 & 170.75
\end{pmatrix}.\]
At first we  obtain  entries of ${\rm \bf A}^{+}$ by (\ref{eq1}):
\[\begin{array}{c}
  d_{3}({\rm \bf A}^{*}{\rm \bf A})=\left|\begin{array}{ccc}
  63 & -44 & -40 \\
  -44 & 48 & -40 \\
  -40 & -40 & 299
\end{array}\right|+\left|\begin{array}{ccc}
  63 & -44 & -11.5 \\
  -44 & 48 & 62 \\
  -11.5 & 62 & 170.75
\end{array}\right|+ \\
  +\left|\begin{array}{ccc}
  63 & -40 & -11.5 \\
  -40 & 299 & -205 \\
  -11.5 & -205 & 170.75
\end{array}\right|+\left|\begin{array}{ccc}
  48 & -40 & 62 \\
  -40 & 299 & -205 \\
  62 & -205 & 170.75
\end{array}\right|=102060,
\end{array}
\]
\[\begin{array}{l}
   l_{11}=\left|\begin{array}{ccc}
  2 & -44 & -40 \\
  0 & 48 & -40 \\
  -5 & -40 & 299
\end{array}\right|+\left|\begin{array}{ccc}
  2 & -44 & -11.5 \\
  0 & 48 & 62 \\
  4 & 62 & 170.75
\end{array}\right|+\left|\begin{array}{ccc}
  2 & -40 & -11.5 \\
  -5 & 299 & -205 \\
  4 & -205 & 170.75
\end{array}\right|=\\
  =25779,
\end{array}
\] and so forth.
Continuing in the same way, we get
\[
{\rm \bf A}^{+}=\frac{1}{102060}\begin{pmatrix}
  25779 & -4905 & 20742 & -5037 \\
  -3840 & -2880 & -4800 & -960 \\
  28350 & -17010 & 22680 & -5670 \\
  39558 & -18810 & 26484 & -13074
\end{pmatrix}.\]
Now we  obtain the least squares solution of the system
(\ref{eq17}) by the matrix method.
 \[{\rm \bf x}^{0}=
   \begin{pmatrix}
  x^0_{11} \\
  x^0_{21} \\
  x^0_{31} \\
  x^0_{41}
\end{pmatrix}=\frac{1}{102060}\begin{pmatrix}
  25779 & -4905 & 20742 & -5037 \\
  -3840 & -2880 & -4800 & -960 \\
  28350 & -17010 & 22680 & -5670 \\
  39558 & -18810 & 26484 & -13074
\end{pmatrix}\cdot\begin{pmatrix}
 1\\
  2 \\
  3 \\
  1
\end{pmatrix}=\]

 \[ =\frac{1}{102060}\begin{pmatrix}
  73158 \\
  -24960 \\
  56700 \\
  68316
\end{pmatrix}=\begin{pmatrix}
  \frac{12193}{17010}\\
  -\frac{416}{1071} \\
 \frac{5}{9} \\
  \frac{5693}{8505}
\end{pmatrix}.
\]
Next we get the least squares solution of the system (\ref{eq17})
by the Cramer rule (\ref{eq14}), where
\[{\rm \bf f}=\begin{pmatrix}
  2 & 7 & 3 & 1 \\
  0 & -4 & -4 & -4 \\
  -5 & -9 & 7 & 12 \\
  4 & 1.5 & -6.5 & -10.5
\end{pmatrix}\cdot\begin{pmatrix}
 1\\
  2 \\
  3 \\
  1
\end{pmatrix}=\begin{pmatrix}
 26\\
  -24 \\
  10 \\
  -23
\end{pmatrix}.
\]
Thus we have
\[
  x^0_1=\frac{1}{102060}\left(
  \left|\begin{array}{ccc}
  26 & -44 & -40 \\
  -24 & 48 & -40 \\
  10 & -40 & 299
\end{array}\right|+
\left|\begin{array}{ccc}
  26 & -44 & -11.5 \\
  -24 & 48 & 62 \\
  -23 & 62 & 170.75
\end{array}\right|+\right.\]
\[\left.+\left|\begin{array}{ccc}
  26 & -40 & -11.5 \\
  10 & 299 & -205 \\
  23 & -205 & 170.75
\end{array}\right|\right)=
\frac{73158}{102060}=\frac{12193}{17010};\]
\[  x^0_2=\frac{1}{102060}\left(\left|\begin{array}{ccc}
  63 & 26 & -40 \\
  -44 & -24 & -40 \\
  -40 & 10 & 299
\end{array}\right|+\left|\begin{array}{ccc}
  63 & 26 & -11.5 \\
  -44 & -24 & 62 \\
  -11.5 & -23 & 170.75
\end{array}\right|+\right. \]\[
 \left.+\left|\begin{array}{ccc}
  -24 & -40 & 62 \\
  10 & 299 & -205 \\
  -23 & -205 & 170.75
\end{array}\right|\right)=\frac{-24960}{102060}=-\frac{416}{1071};
\]
\[
  x^0_3=\frac{1}{102060}\left(\left|\begin{array}{ccc}
  63 & -44 & 26 \\
  -44 & 48 & -24 \\
  -40 & -40 & 10
\end{array}\right|
  +\left|\begin{array}{ccc}
  63 & 26 & -11.5 \\
  -40 & 10 & -205 \\
  -11.5 & -23 & 170.75
\end{array}\right|+\right.\]\[\left.+\left|\begin{array}{ccc}
  48 & -24 & 62 \\
  -40 & 10 & -205 \\
  62 & -23 & 170.75
\end{array}\right|\right)=\frac{56700}{102060}=\frac{5}{9};
\]
\[
  x^0_4=\frac{1}{102060}\left(\left|\begin{array}{ccc}
  63 & -44 & 26 \\
  -44 & 48 & -24 \\
  -11.5 & 62 & -23
\end{array}\right| \\
  +\left|\begin{array}{ccc}
  63 & -40 & 26 \\
  -40 & 299 & 10 \\
  -11.5 & -205 & -23
\end{array}\right|+\right.\]
\[\left.+\left|\begin{array}{ccc}
  48 & -40 & -24 \\
  -40 & 299 & 10 \\
  62 & -205 &-23
\end{array}\right|\right)=\frac{68316}{102060}=\frac{5693}{8505}.
\]

2. Let us consider the following system of linear equations.
\begin{equation} \label{eq18}
\left\{ \begin{array}{c}
x_1-x_2+x_3+x_4=1, \\
  x_2-x_3+x_4=2, \\
  x_1-x_2+x_3+2x_4=3, \\
  x_1-x_2+x_3+x_4=1.
\end{array}\right.\end{equation}
The coefficient matrix of the system is the matrix ${\rm \bf
A}=\begin{pmatrix}
  1 & -1 & 1 & 1 \\
  0 & 1 & -1 & 1 \\
  1 & -1 & 1 & 2 \\
  1 & -1 & 1 & 1
\end{pmatrix}$. It is easy to verify the following:
\[{\rm \bf
A}^{2}=\begin{pmatrix}
  3 & -4 & 4 & 3 \\
  0 & 1 & -1 & 0 \\
  4 & -5 & 5 & 4 \\
  3 & -4 & 4 & 3
\end{pmatrix},\,\,{\rm \bf
A}^{3}=\begin{pmatrix}
  10 & -14 & 14 & 10 \\
  -1 & 2 & -2 & -1 \\
  13 & -18 & 18 & 13 \\
  10 & -14 & 14 & 10
\end{pmatrix},\] and $\rank{\rm \bf A}=3$,
 $\rank{\rm {\bf A}}^{2}=\rank{\rm \bf A}^{3}=2$.
 This implies $k=Ind({\rm {\bf A}})=2$.
We  obtain  entries of ${\rm \bf A}^{D}$ by   (\ref{eq11}).
\[\begin{array}{c}
   d_{2}({\rm \bf A}^{3})=\left|\begin{array}{cc}
  10 & -14 \\
  -1 & 2
\end{array}\right|+\left|\begin{array}{cc}
  10 & 14 \\
  13 & 18
\end{array}\right|+\left|\begin{array}{cc}
  10 & 10 \\
  10 & 10
\end{array}\right|\\
   +\left|\begin{array}{cc}
 2 & -2 \\
  -18 & 18
\end{array}\right|+\left|\begin{array}{cc}
  2 & -1 \\
  -14 & 10
\end{array}\right|+\left|\begin{array}{cc}
  18 & 13 \\
  14 & 10
\end{array}\right|=8,
\end{array}
\]
\[   d_{11}=\left|\begin{array}{cc}
 3 & -14 \\
  0 & 2
\end{array}\right|+\left|\begin{array}{cc}
  3 & 14 \\
  4 & 18
\end{array}\right|+\left|\begin{array}{cc}
  3 & 10 \\
  3 & 10
\end{array}\right|=4,
\] and so forth.
Continuing in the same way, we get ${\rm \bf
A}^{D}=\begin{pmatrix}
  0.5 & 0.5 & -0.5 & 0.5 \\
 1.75 & 2.5 & -2.5 & 1.75 \\
  1.25 & 1.5 & -1.5 & 1.25 \\
  0.5 & 0.5 & -0.5 & 0.5
\end{pmatrix}.$
Now we obtain the Drazin inverse solution ${\rm {\bf
\widehat{x}}}$ of the system (\ref{eq18}) by the Cramer rule
(\ref{eq16}), where
\[{\rm \bf g}={\rm \bf
A}^{2}{\rm {\bf y}}=\begin{pmatrix}
  3 & -4 & 4 & 3 \\
  0 & 1 & -1 & 0 \\
  4 & -5 & 5 & 4 \\
  3 & -4 & 4 & 3
\end{pmatrix}\cdot\begin{pmatrix}
 1\\
  2 \\
  3 \\
  1
\end{pmatrix}=\begin{pmatrix}
 10\\
  -1 \\
  13 \\
  10
\end{pmatrix}.
\]
Thus we have
\[
  \widehat{x}_1=\frac{1}{8}\left(
  \left|\begin{array}{cc}
 10 & -14 \\
  -1 & 2
\end{array}\right|+\left|\begin{array}{cc}
  10 & 14 \\
  13 & 18
\end{array}\right|+\left|\begin{array}{cc}
  10 & 10 \\
  10 & 10
\end{array}\right|\right)=\frac{1}{2},
\]
\[
  \widehat{x}_2=\frac{1}{8}\left(
  \left|\begin{array}{cc}
 10 &10 \\
  -1 & -1
\end{array}\right|+\left|\begin{array}{cc}
  -1 & -2 \\
  13 & 18
\end{array}\right|+\left|\begin{array}{cc}
 -1 & -1 \\
  10 & 10
\end{array}\right|\right)=1,
\]
\[
  \widehat{x}_3=\frac{1}{8}\left(
  \left|\begin{array}{cc}
 10 &10 \\
  13 & 13
\end{array}\right|+\left|\begin{array}{cc}
  2 & -1 \\
  -18 & 13
\end{array}\right|+\left|\begin{array}{cc}
13 & 13 \\
  10 & 10
\end{array}\right|\right)=1,
\]
\[
  \widehat{x}_4=\frac{1}{8}\left(
  \left|\begin{array}{cc}
 10 & 10 \\
  10 & 10
\end{array}\right|+\left|\begin{array}{cc}
  2 & -1 \\
  -14 & 10
\end{array}\right|+\left|\begin{array}{cc}
  18 & 13 \\
  14 & 10
\end{array}\right|\right)=\frac{1}{2}.
\]

{\bf Acknowledgment.} The author would like to thank Professor
Chi-Kwong Li and the referee for their useful suggestions.

\end{document}